\documentclass[a4paper,12pt]{article}
\usepackage{amssymb}
\usepackage{amsthm}
\usepackage{tikz}
\usepackage{tikz,pgfplots}
\usetikzlibrary{decorations.markings}
\usepackage{amsmath,amssymb}
\usepackage{todonotes}
\usepackage{url}
\usepackage{mathtools,booktabs}

\usepackage{color}
\usepackage{graphicx}
\usepackage{placeins,caption} 

\usepackage{hyperref}

\hypersetup{colorlinks=true, linkcolor=blue, citecolor=blue, urlcolor=blue}

\DeclareMathOperator{\mob}{Mob_{gp}}

\DeclareMathOperator{\gp}{gp}

\DeclareMathOperator{\cp}{\,\square\,}

\DeclareMathOperator{\cmob}{Mob_{gp}^*}

\DeclareMathOperator{\rob}{\mathsf{R}}

\usepackage{stackengine}

\usepackage[top=3cm,bottom=3cm,left=2.8cm,right=2.8cm]{geometry}
	\newtheorem{theorem}{Theorem}[section]
	\newtheorem{lemma}[theorem]{Lemma}
	
	\newtheorem{proposition}[theorem]{Proposition}

	\theoremstyle{definition}
	
	\newtheorem{claim}[theorem]{Claim}
	\newtheorem{definition}[theorem]{Definition}

	\newcommand{\move}{\rightsquigarrow}

 \pgfplotsset{compat=1.18}
	
\begin{document}
	
		\title{Mobile general position for torus and cylinder graphs}
		\author{Ethan Shallcross $^{a}$ \\ \texttt{\footnotesize ethan.shallcross@open.ac.uk} \\ 
		{\footnotesize ORCID: 
		0009-0000-5730-6834}
  \and 
		James Tuite $^{a,b,c}$ \\ \texttt{\footnotesize james.t.tuite@open.ac.uk} \\ 
		{\footnotesize ORCID: 
		0000-0003-2604-7491}
  \and 
        Aleksander Vesel $^{d,e}$ \\ \texttt{\footnotesize aleksander.vesel@um.si} \\ 
		{\footnotesize ORCID: 
		0000-0003-3705-0071}
}

\maketitle

\begin{center}
$^a$ School of Mathematics and Statistics, Open University, Milton Keynes, UK\\

$^b$ Department of Informatics and Statistics, Klaip\.{e}da University, Lithuania\\

$^c$ Faculty of Mathematics and Informatics, Vilnius University, Lithuania\\

$^d$ Faculty of Natural Sciences and Mathematics,
University of Maribor, Slovenia  \\

$^e$ IMFM, Ljubljana, Slovenia

\end{center}

\begin{abstract}
A set of vertices of a graph is in general position if no three of them lie on a common shortest path. In the mobile general position problem, we imagine a collection of robots on the vertices of a graph that must move along the edges one at a time, always remaining in general position, so that every vertex is visited by a robot. The mobile general position number of a graph is the largest number of robots for which this is possible. Following earlier investigations of this property in Cartesian products, this paper solves the open problem of determining its value for all cylinder and torus graphs. We also consider the completely mobile variant in cylinders, in which every robot must visit every vertex. 
\end{abstract}


\noindent
AMS Subj.\ Class.\ (2020): 05C12, 05C69, 05C76

\section{Introduction}\label{sec:intro}

The general position problem for graphs, introduced in~\cite{ullas-2016} and~\cite{manuel-2018}, asks for a largest vertex subset $S$ of a graph $G$ such that no three vertices of $S$ lie on a common shortest path of $G$. A survey of the problem is given in~\cite{survey}. A subset with this property is a \emph{general position set} of $G$ and the number of vertices in a largest general position set of $G$ is the \emph{general position number} $\gp (G)$ of $G$. A useful visualisation of the problem was suggested in the arXiv version of~\cite{manuel-2018}: robots are stationed at the vertices of a graph, with at most one robot per vertex, and they communicate with each other by sending signals along shortest paths. To avoid any communication channels being blocked, we require that no third robot lies on a shortest path between a pair of other robots.

A dynamic version of this problem was introduced in~\cite{klavzar-2023}. In the so-called \emph{mobile general position problem} a swarm of robots has to move through the graph, with just one robot moving at a time, in such a way that the robots are always in general position and every vertex is visited by a robot at some time. The largest number of robots in such a swarm is the \emph{mobile general position number} of $G$, denoted $\mob (G)$. A more restrictive version, called the \emph{completely mobile general position problem}, was presented in~\cite{Boshar2024}, in which each robot is required to visit every vertex, and the largest possible number of robots in such a swarm is called the \emph{completely mobile general position number} of $G$, written $\cmob (G)$. Mobile versions of well-known graph theory problems have since been researched, such as mobile mutual visibility~\cite{Dettlaf2024} and mobile metric resolving sets~\cite{Camacho,CamachoMoreno}. A similar independent direction involving dynamic $d$-independent sets can be found in~\cite{span}.

The general position problem has been investigated extensively for Cartesian product graphs. Some examples are~\cite{klavzar-2021,KorzeVesel,Manuel-2018b,tian-2021a,tian-2021b}. The mobile general position problem for Cartesian products was addressed in~\cite{klavzar-2025}. As well as giving some general bounds for the mobile general position number of a Cartesian product $G \cp H$, the article solved the problem for the Cartesian product of a complete graph with a path, the product of two complete graphs (i.e.\ a rook graph), and the product of a tree with a complete graph $K_2$ with order two. The paper also solved the problem for most cylinder graphs $C_r \cp P_s$:

\begin{proposition}[\cite{klavzar-2025}]
The mobile general position number of $C_r \cp P_s$ is given by
\[\mob(C_r \cp P_s) = \begin{cases}
    2, & \text{if}\ (r,s) = (4,2),  \\
    3, & \text{if}\ r = 4 \text{ and } s \geq 3, \text{ or } r = 3 \text{ and } s \geq 2, \\
    4, & \text{if}\ r \in \{ 5,6,8\} \text{ and } s \geq 2, \text{ or } r \geq 5 \text{ and } s \leq 4, \\
    5, & \text{if}\ s \geq 5 \text{ and either } r = 9 \text{ or } r \geq 11.
\end{cases} \]
\label{prop:knowncylinders}
\end{proposition}
Hence the only open cases for cylinder graphs were $C_7 \cp P_s$ and $C_{10} \cp P_s$ for $s \geq 5$, and the article conjectured that the answer is four in these cases. Much less was done in~\cite{klavzar-2025} on torus graphs $C_r \cp C_s$, the only significant result being the following lower bound.

\begin{proposition}[\cite{klavzar-2025}]
    For $r=9$ or $r > 10$, and $s \geq 10$, $\mob(C_r \cp C_s) \geq 5$.
\end{proposition}

It is also stated in~\cite{klavzar-2025} that $\mob(C_9 \cp C_8) = 7$. In this article we settle the remaining cases. Section~\ref{sec:cylinders} begins by finding the mobile general position numbers of the cylinders $C_7 \cp P_s$ and $C_{10} \cp P_s$. We then determine the completely mobile general position numbers of all cylinder graphs. Section~\ref{sec:twocycle} begins by finding the mobile general position numbers of all sufficiently large torus graphs. We then complete the classification of mobile general position numbers for graphs in this family, resolving another question from~\cite{klavzar-2025}.

For basic graph theory terminology we refer to~\cite{graphsanddigraphs}. A few words on the notation that we will use in the article. For positive integers $n$ and $k$, $n \not = k$, we will use the notation $[k] = \{1, 2, \ldots, k\}$, 
$[k]_0 = \{0, 1, \ldots, k-1\}$ and 
\begin{displaymath}
[k,n] =
        \left \{ \begin{array}{llll}
              \{k, k+1, \ldots, n\}, &  n \ge k   \\
              \{n, n+1, \ldots, k\}, &  k > n    \\
             \end{array}. \right.
\end{displaymath}

For a positive integer $s$ and $n,k \in [s]_0$  such that $n,k < s $ and $n \not = k$,
  let $$|n - k|_s = \min \{ |n - k|, s - |n - k|\}.$$ 
Moreover, let
\begin{displaymath}
[k,n]_s =
        \left \{ \begin{array}{llll}
              [k,n], & \text{ if } |n - k| <  s - |n - k|     \\
              {[s]_0}, & \text{ if } |n - k| =  s - |n - k|   \\
              ({[s]_0} \setminus [k,n]) \cup \{n,k\}, &  \text{ if } \;  |n - k| >  s - |n - k|   \\
             \end{array}. \right.
\end{displaymath}

The  \emph{Cartesian product}  $G \cp H$ of graphs $G$ and $H$ has the vertex set $V(G \cp H)=V(G)\times V(H)$,  where vertices $(g,h)$ and $(g',h')$ are adjacent if either $gg' \in E(G)$ and $h=h'$, or $g=g'$ and $hh' \in E(H)$. Given a vertex $h \in V(H)$, the subgraph of $G \cp H$ induced by the set $ \{(g,h) : g \in V(G)\}$ is called a $G$-\emph{layer} and is denoted by $G^h$. Analogously, $H$-\emph{layers} ${}^{g}H$ are defined for vertices $g \in V(G)$.

We will assume that $V(C_s)=[s]_0$ for every  $s\ge 3$. Thus,  
if $u \in V(C_s \cp C_t)$,  we write $u=(u_x,u_y)$, where $u_x \in [s]_0$ and  $u_y \in [t]_0$. 

The following lemma is a direct consequence of \cite[Lemma 3.1(ii),(iii)]{klavzar-2021}
\begin{lemma}
\label{lem:2-in-a-layer}
Let $r \ge 2$, $s \ge 3$, and let $S$ be a general position set of $P_r \cp C_s$. If $|S \cap {}^{i}C_s| \ge 2$ for some $i \in [r]_0$, then $|S| \le 4$.
\end{lemma}

\begin{proposition}
\label{prop:one-per-layer}
Let $r,s \ge 7$. If $S$ is a general position set of $C_r \cp C_s$ with $|S| = 7$, then every $C_r$-layer and every $C_s$-layer of $C_r \cp C_s$ contains at most one vertex of $S$.
\end{proposition}

\begin{proof}
By the commutativity of the Cartesian product of graphs, it suffices to prove the assertion for $C_s$-layers. Suppose, to the contrary, that some $C_s$-layer contains at least two vertices of $S$. As $C_r$ is vertex-transitive, we may assume that this layer is ${}^{0}C_s$.

Let $k = \lfloor r/2 \rfloor$ and 
$A = \{0,1,\ldots,k\}$, $B = \{k+1,k+2,\ldots,r-1,0\}$.
Then $A \cup B = [r]_0$ and $A \cap B = \{0\}$. In $C_r$, the set $A$ induces a path on $k+1$ vertices, while $B$ induces the path $k+1,k+2,\ldots,r-1,0$ on $r-k = \lceil r/2 \rceil$ vertices. Each of these two paths has at most $\lfloor r/2 \rfloor$ edges, so it is a geodesic of $C_r$ and hence an isometric subgraph of $C_r$. Thus, the subgraphs $G_A$ and $G_B$ of $C_r \cp C_s$ induced by $A \times [s]_0$ and $B \times [s]_0$, respectively, are isometric subgraphs of $C_r \cp C_s$. Moreover, $G_A$ and $G_B$ are 
isomorphic to  $P_{k+1} \cp C_s$ and  $P_{r-k} \cp C_s$, respectively.

Since every geodesic of an isometric subgraph is also a geodesic of $C_r \cp C_s$, the sets $S_A = S \cap (A \times [s]_0)$ and $S_B = S \cap (B \times [s]_0)$ are general position sets of $G_A$ and $G_B$, respectively. The set ${}^{0}C_s$ is a $C_s$-layer of both $G_A$ and $G_B$, and it contains at least two vertices of $S_A$ as well as of $S_B$. Hence, Lemma~\ref{lem:2-in-a-layer} yields $|S_A| \le 4$ and $|S_B| \le 4$. Finally, $S = S_A \cup S_B$ and $S_A \cap S_B = S \cap {}^{0}C_s$, and therefore
\[
|S| = |S_A| + |S_B| - |S \cap {}^{0}C_s| \le 4 + 4 - 2 = 6,
\]
a contradiction with $|S| = 7$.
\end{proof}

\section{Mobile general position for cylinder graphs}\label{sec:cylinders}

In this section we examine the mobile and completely mobile general position problems for cylinder graphs. In Subsection~\ref{subsec:mob of cylinder} we completely solve the remaining open cases for the mobile general position problem on cylinders and in Subsection~\ref{subsec:comp mob cylinders} we consider completely mobile general position numbers of cylinders.

\subsection{Mobile general position numbers}\label{subsec:mob of cylinder}

We first find $\mob(C_r \cp P_s)$ for $r = 7, 10$ and $s \geq 5$, settling an open problem from~\cite{klavzar-2025}. Theorem 3.4 and Lemma 3.5 of~\cite{klavzar-2025} together imply that the mobile general position numbers of these graphs are at least 4, and Theorem 3.2 of~\cite{klavzar-2021} shows that the largest possible general position number of a cylinder is 5, so it remains only to rule out the existence of a mobile general position set of size 5 in $C_r \cp P_s$ for $r = 7,10$ and $s \geq 5$. 

\begin{lemma}\label{lem:distinct coordinates}
If $S$ is a general position set with $|S| \geq 5$ of $C_r \cp P_s$, where $r \in \{ 7,10\} $ and $s \geq 5$, then all first coordinates of vertices in $S$ are distinct, as are all second coordinates.
\end{lemma}
\begin{proof}
Suppose that $S$ is a general position set of $C_r \cp P_s$ of order $|S| = 5$, where $r \in \{ 7,10\} $ and $s \geq 5$. By Lemma~\ref{lem:2-in-a-layer}
there must be at most one vertex of $S$ in every $C_r$-layer. 

We claim that all of the first coordinates of $S$ must be distinct, i.e. no two robots lie in the same $P_s$-layer. Suppose that $S$ has two elements with the same first coordinate. If $S \cap ([r]_0 \times \{ j \}) = \emptyset$ and $(x,j+1) \in S$, then the move $(x,j+1) \move (x,j)$ is legal. In this way, by shifting robots through empty layers to smaller coordinates, $S$ can be transformed into a general position set of $C_r \cp P_5$ with the same property. Hence it suffices to prove that there is no such general position set in $C_r \cp P_5$. We can thus write $S$ as $\{ (x_i,i):0 \leq i \leq 4\} $. 

Suppose that there are two vertices of $S$ that lie in the same $P_5$-layer, say $(x,i)$ and $(x,j)$, where $i < j$. If $j < 4$, then there would be a shortest path from the robot at $(x,i)$ to the robot in $[r]_0 \times \{ 4\} $ via the robot at $(x,j)$. Similarly we cannot have $i > 0$. Hence we can assume that $(0,0),(0,4) \in S$. By the preceding argument the remaining three robots in $[r]_0 \times [1,3]$ must have distinct first coordinates. Set $A(7) = [1,3] $, $B(7) = [4,6] $, $A(10) = [1,5] $, $B(10) = [5,9]$. As there are three more robots, either $A(r)$ or $B(r)$ must contain the first coordinates $x_i,x_j$ of at least two of the robots. By rotational and reflectional symmetry of the cylinder, we can assume that these coordinates lie in $A(r)$, $x_i < x_j$ and $i < j$. But then there is a shortest path from $(x_j,j)$ to $(0,0)$ via $(x_i,i)$, a contradiction.
\end{proof}

\begin{lemma}[Scrunching Lemma]
For $r \in \{ 7,10\} $ and $s \geq 5$, $\mob (C_r \cp P_s) = \mob (C_r \cp P_5)$.
\end{lemma}
\begin{proof}
 If $\mob (C_r \cp P_s) = 4$ for $r \in \{ 7,10\} $ and all $s \geq 5$, then our claim holds. Suppose that there is some $s \geq 5$ such that $(C_r \cp P_s)$ has a mobile general position set of size 5 for some $r \in \{ 7,10\} $. By Lemma~\ref{lem:distinct coordinates} all the first coordinates and all the second coordinates are distinct. Hence we can write $S = \{ (x_i,y_i) : 0 \leq i \leq 4\}$, where $x_0 < x_1 < \ldots < x_4$. Now observe that the general position property for $S$ does not depend on the exact values of the second coordinates, but only on the ordering of the sequence $y_0,y_1,y_2,y_3,y_4$.  

 For any configuration $S$ of the five robots in $C_r \cp P_s$, define the `scrunched' version $S'$ of $S$ to be the result of performing moves of the form $(i,j) \move (i,j-1)$ whenever the $(j-1)$-layer is empty, until no further such moves can be performed. By the above observation, $S'$ is in general position if and only if $S$ is in general position. Suppose that $S_0,S_1,\ldots $ is a sequence of configurations of the robots in $C_r \cp P_s$ linked by legal moves that visits every vertex of $C_r \cp P_s$ and replace each $S_i$ by its scrunched version $S_i'$. Then remove any consecutive pairs for which $S_i' = S_{i+1}'$. The resulting sequence will allow the robots to visit every vertex of $C_r \cp P_5$ in general position. Thus $\mob (C_r \cp P_s) = \mob (C_r \cp P_5)$.
\end{proof}

\begin{proposition}
    For all $s \geq 5$, $\mob(C_7 \cp P_s) = \mob(C_{10} \cp P_s) = 4$.
\end{proposition}
\begin{proof}
By the preceding results, it suffices to show that $\mob (C_r \cp P_5) < 5$ for $r \in \{ 7,10\} $. Suppose that $S$ is a mobile general position set in $C_r \cp P_5$ with five robots. By Lemma~\ref{lem:distinct coordinates} each robot can only move in its $C_r$-layer. We focus on the configurations before and after the move $(0,0) \move (1,0)$. 

Consider the configuration of robots when a robot is at $(0,0)$. Let $A = [1,\left \lfloor \frac{r}{2} \right \rfloor ] \times [5]_0$ and $B = [\left \lceil \frac{r}{2} \right \rceil ,r-1] \times [5]_0$. Suppose that one of these sets, say $A$, contains three robots. Let their positions be $(i_1,j_1),(i_2,j_2),(i_3,j_3)$, where $i_1 < i_2 < i_3$. To avoid creating three-in-line with $(0,0)$ we must have $j_1 > j_2 > j_3$; however, this implies that $(i_2,j_2)$ lies on a shortest path from $(i_1,j_1)$ to $(i_3,j_3)$. The same argument applies to $B$, so $A$ and $B$ contain two robots each. When $r = 10$, the overlap of $A$ and $B$ at $5$ implies that there is no robot with first coordinate $5$.

Applying the same argument to the moment after the move, when there is a robot at $(1,0)$, we see that when $r = 7$ the remaining four robots have coordinates $(2,j_1),(3,j_2),(5,j_3),(6,j_4)$ and when $r = 10$ we can write the positions of the robots as  \[ (i_1,j_1),(i_2,j_2),(i_3,j_3),(i_4,j_4),\] where $i_1,i_2 \in [2,4]$, $i_3,i_4 \in [7,9]$ and $i_1<i_2<i_3<i_4$. In both cases to avoid three-in-line with $(1,0)$ we must have $j_1 > j_2$ and $j_3 < j_4$.

For $r = 7$, to avoid $(3,j_2)$ lying on a shortest path from $(2,j_1)$ to $(5,j_3)$ we must have $j_2 < j_3$, but this implies that there is a shortest path from $(3,j_2)$ to $(6,j_4)$ via $(5,j_3)$, a contradiction. Now assume that $r = 10$. Suppose that $j_2 > j_3$, implying that $j_1 > j_2 > j_3$. Then to avoid a shortest $(i_1,j_1),(i_3,j_3)$-path through $(i_2,j_2)$ we must have $i_1 = 2, i_3 = 8, i_4 = 9$. To avoid a shortest $(i_3,j_3),(i_1,j_1)$-path through $(i_4,j_4)$, it follows that $j_1 < j_4$. But now there is a shortest $(i_4,j_4),(i_2,j_2)$-path through $(i_1,j_1)$. Hence $j_2 < j_3$. Then to avoid a shortest $(i_2,j_2),(i_4,j_4)$-path through $(i_3,j_3)$ we have $i_1 = 2, i_2 = 3, i_4 = 9$. To avoid a shortest $(i_4,j_4),(i_2,j_2)$-path through $(i_1,j_1)$, we need $j_1 > j_4$; but this implies that there is a shortest $(i_3,j_3),(i_1,j_1)$-path through $(i_4,j_4)$. This contradiction concludes our proof.
\end{proof}

\subsection{Completely mobile general position numbers}\label{subsec:comp mob cylinders}

We now extend our investigation into mobile general position problems in cylinder graphs to consider the completely mobile variant. In this subsection, we shall completely classify the completely mobile general position numbers of all cylinder graphs.

We begin with the cylinder graph $C_r \cp K_2$, also known as the prism graphs. Recall from Proposition~\ref{prop:knowncylinders} that \[
\mob(C_r \cp K_2) = \begin{cases}
3, & \text{if}\ r = 3,  \\
2, & \text{if}\ r = 4,  \\
4, & \text{otherwise}.  \\
\end{cases} \]
We now find the completely mobile general position numbers of these graphs.

\begin{proposition}
For $r \geq 3$,
\[\cmob(C_r \cp K_2) = \begin{cases}
    2, & \text{if}\ r \leq 4 \\
    3, & \text{if}\ r > 4.
\end{cases} \]
\label{prop:cmobsmallcylinders}
\end{proposition}
\begin{proof}
By~\cite[Lemma 3.1]{klavzar-2021} any general position set of $C_r \cp K_2$ of size at least $4$ has $\leq 2$ vertices in each $C_r$-layer. Therefore, if $4$ robots are stationed on $C_r \cp K_2$, no robot may move to the other $C_r$-layer and the robots fail to be completely mobile. Hence $\cmob(C_r\cp K_2) \leq 3$ for $r \geq 3$. It is shown in~\cite{klavzar-2023} that $\mob (C_r) = 3$ when $r \geq 7$ or $r \in \{ 3,5\} $. When $r = 5$ or $r \geq 7$ by iterating the sequence of moves used in the proof of that result, every robot visits every vertex of $C_r$, and so $\cmob (C_r) = 3$ for $r = 5$ and $r \geq 7$. Identify a $C_r$-layer of $C_r \cp K_2$ with $C_r$ and perform this sequence of moves; then it can be seen that the three robots can move into the other $C_r$-layer one at a time without violating the general position property, and then perform the same sequence of moves in the other $C_r$-layer. Thus $\cmob (C_r \cp K_2) = 3$ when $r = 5$ or $r \geq 7$.  

It remains to prove that $\cmob(C_6 \cp K_2) \geq 3$. Station robots $\rob_0, \rob_1, \rob_2$ on $C_6 \cp K_2$ at $(0,0)$, $(3,0)$, $(1,1)$, respectively. By repeatedly performing moves of the form $(i,j) \move (i+1, j)$ on the robots $\rob_2, \rob_1, \rob_0$ in that order, each robot can visit each vertex in its $C_6$-layer whilst maintaining general position. For the robots to swap layers, from the starting configuration perform the sequence of moves $R_2: (1,1) \move (2,1)$, $R_1: (3,0) \move (4,0) \move (4,1)$, $R_0: (0,0) \move (0,1)$, $R_2: (2,1) \move (2,0)$ and $R_1: (4,1) \move (3,1)$. This switches the layer of each robot and results in a configuration isomorphic to the starting configuration, so now each robot can visit each vertex in the other $C_6$-layer. Thus $\cmob (C_6 \cp K_2) = 3$.

It follows immediately from~\cite[Theorem 3.4]{klavzar-2025} that \[ \cmob(C_4 \cp K_2) \leq \mob(C_4 \cp K_2) = 2,\] and so $\cmob(C_4 \cp K_2) = 2$. In $C_3 \cp K_2$ a subset of size $3$ is in general position if and only if all the first coordinates are different, so no robot in such a set is able to move to a vertex with a different first coordinate, and it follows that $\cmob(C_3 \cp K_2) = 2$. 
\end{proof}

The completely mobile general position numbers of the remaining cylinder graphs are classified in the following result. We need a lemma that shows that the completely mobile general position number of the cylinder is monotonic in the length of the path factor.

\begin{lemma}\label{lem:shifting}
For $r,s \geq 3$ and $t \geq 0$, \[ \cmob (C_r \cp P_{s+t}) \geq \cmob (C_r \cp P_s).\] 
\end{lemma}
\begin{proof}
    Suppose that $S$ is a completely mobile general position set in $C_r \cp P_s$ with $|S| = \cmob (C_r \cp P_s)$. We can label the robots $R_0,R_1,\ldots ,R_{|S|-1}$ in non-decreasing order in the second coordinate. The cylinder $C_r \cp P_s$ is isomorphic to the isometric subgraph of $C_r \cp P_{s+t}$ induced by $C_r \times [s]_0$. Station the robots of $S$ at the corresponding positions in $C_r \cp P_{s+t}$, move the robots around the copy of $C_r \cp P_s$ so that all vertices of $[r]_0 \times [s]_0$ are visited in $C_r \cp P_{s+t}$ and then return the robots to the starting configuration. Then perform a move of the form $(i,j) \move (i,j+1)$ for the robots in the order $R_{|S|-1},R_{|S|-2},\ldots ,R_1,R_0$. All of these moves are legal. This moves the robots into an isomorphic configuration in an isometric copy of $C_r \cp P_s$ on the vertices $[r]_0 \times [1,s]$, and the robots can now repeat their move sequence to visit all the vertices of this subgraph. By repeating this argument the robots can repeat their move sequence in any subgraph induced by $[r]_0 \times [p,p+s-1]$ for $0 \leq p \leq t$ and we see that the configuration was completely mobile.
\end{proof}

\begin{theorem}
For all integers $r, s \geq 3$, \[
\cmob(C_r \cp P_s) = \begin{cases}
    3, & \text{if}\ 3 \leq r \leq 4 \\
	4, & \text{if}\ r > 4,
\end{cases}
\]
except that $\cmob(C_6 \cp P_3) = 3$.
\end{theorem}
\begin{proof}
From the classification of the mobile general position numbers of cylinders, we have immediately that the completely mobile general position number of a cylinder can be at most five. If there are five robots in a completely mobile general position set of a cylinder, then we see from Lemma~\ref{lem:2-in-a-layer} 
that there is at most one robot per $C_r$-layer; however, this implies that the robots would not be able to change the ordering of their second coordinates, a contradiction. Hence  $\cmob(C_r\cp P_s) \leq 4$ for $r,s \geq 3$.

We start with small $r$. It follows from Propositions 2.2 and 3.6 of~\cite{klavzar-2025} that the completely mobile general position numbers of $C_3 \cp P_s$ and $C_4 \cp P_s$ are at most $3$. By Lemma~\ref{lem:shifting} it suffices to show that $\cmob (C_3 \cp P_3) \geq 3$ and $\cmob (C_4 \cp P_3) \geq 3$. For $r = 3$ the set $\{ (0,0),(1,1),(2,2)\}$ is easily seen to be a mobile general position set, and the positions of the robots can be permuted cyclically through the moves $(0,0) \move (0,1)$, $(1,1) \move (1,0)$, $(2,2) \move (1,2)$, $(0,1) \move (2,1)$, $(1,0) \move (0,0)$, $(1,2) \move (1,1)$ and $(2,1) \move (2,2)$. For $r = 4$ start with the configuration $\{ (0,0),(2,1),(0,2)\} $. By shifting the robots around the cycle layers, it can be seen that the robots are a mobile general position set, and to permute the positions of the robots, perform the moves $(0,2) \move (3,2)$, $(2,1) \move (2,0)$, $(0,0) \move (0,1)$, $(3,2) \move (2,2)$, $(2,0) \move (1,0)$, $(0,1) \move (0,2)$, $(2,2) \move (2,1)$, $(1,0) \move (0,0)$, so that they are completely mobile.

Now we turn to the exceptional cylinder $C_6 \cp P_3$. By applying Lemma~\ref{lem:shifting} to Proposition~\ref{prop:cmobsmallcylinders} we see that $\cmob(C_6 \cp P_3) \geq 3$. Suppose for a contradiction that there is a completely mobile general position set in this cylinder with four robots. Recall from~\cite{klavzar-2021} that there are always at most two robots per $C_6$-layer. Also, as $C_6 \cp P_3$ is bipartite, it follows from the results of~\cite{lower,lowercartesian} that the robots must always occupy an independent set. In order for a robot $\rob $ to move from the layer $[6]_0 \times \{ 0\} $ to the layer $[6]_0 \times \{ 2\} $, at some point there must be one robot in each of $[6]_0 \times \{ 0\} $ and $[6]_0 \times \{ 2\} $ and two robots, one of which is $\rob $, in $[6]_0 \times \{ 1\} $. Without loss of generality, suppose that the move of $\rob $ from the middle layer to $[6]_0 \times \{ 2\} $ is $(0,1) \move (0,2)$. Suppose that the other robot in the middle layer is at distance two from $(0,1)$, say it is at $(2,1)$. From the two restrictions mentioned above applied to the configurations before and after this move, the only position that the other robot in $[6]_0 \times \{ 2\} $ can occupy is $(4,2)$ and the robot in $[6]_0 \times \{ 0\} $ is either at $(4,0)$ or $(1,0)$. However, there is a shortest $(4,0),(0,2)$-path through $(4,2)$ and a shortest $(1,0),(4,2)$-path through $(2,1)$, so such a move is impossible. Finally suppose that when $R$ makes the move $(0,1) \move (0,2)$ the other robot in the middle layer is at distance three from $(0,1)$, i.e.\ is at $(3,1)$. By the above restrictions we can suppose that the robot in $[6]_0 \times \{ 2\} $ is at $(2,2)$. However, this would create a shortest $(0,2),(3,1)$-path via $(2,2)$. Thus $\cmob (C_6 \cp P_3) = 3$.

Now we complete the case $r = 6$. Consider $C_6 \cp P_4$. Begin by stationing robots $\rob_0, \dots, \rob_3$ at $(0,0)$, $(2,1)$, $(4,2)$, and $(0,3)$, respectively. By repeatedly moving $\rob_3, \rob_2, \rob_1, \rob_0$ in that order with moves $(i,j) \move (i+1, j)$ around the cycle layers, we see that this is a mobile general position set. To permute the positions of $\rob_0, \rob_1, \rob_2$ in a 3-cycle, move $\rob_2$ by $(4,2) \move (4,1)$, $\rob_1$ by $(2,1) \move (2,2)$, $\rob_0$ by $(0,0) \move (1,0)$, $\rob_2$ by $(4,1) \move (4,0)$, $\rob_0$ by $(1,0) \move (1,1)$, $\rob_1$ by $(2,2) \move (3,2)$, $\rob_0$ by $(1,1) \move (2,1)$, $\rob_2$ by $(4,0) \move (5,0)$, $\rob_1$ by $(3,2) \move (4,2)$, and finally $\rob_2$ by $(5,0) \move (0,0)$. By symmetry (specifically the automorphism $(i,j) \rightarrow (-i,3-j)$), we may also permute the positions of $\rob_1, \rob_2, \rob_3$ in a 3-cycle, and this is sufficient for any robot to occupy any of the positions in the starting configuration, giving $\cmob(C_6 \cp P_4) = 4$. The equality for $s > 4$ follows from Lemma~\ref{lem:shifting}. 

We finally consider the case $r = 5$ or $r \geq 7$. We start by showing $\cmob (C_r \cp P_3) \geq 4$. Set $\alpha = \left \lfloor \frac{r-1}{2} \right \rfloor $. Station robots $\rob_0, \dots, \rob_3$ at $(r-\alpha ,1)$, $(0,0)$, $(0,2)$, $(\alpha ,1)$ respectively. This set is in general position. Perform the moves $(i,j) \move (i+1,j)$ on the robots $\rob _0$, $\rob _1$, $\rob _2$, $\rob _3$ in that order. All of these moves are legal, and by repeating this sequence all vertices of $C_r \cp P_3$ are visited, so the starting configuration is a mobile general position set. 

Now we show that the starting configuration is completely mobile. For $i = \alpha ,\alpha -1 ,\ldots ,3$ in that order perform the pair of moves $\rob _3: (i,1) \move (i-1,1), \rob _1: (r-\alpha +i,0) \move (r-\alpha +i-1,0)$. This brings robot $\rob _3$ to $(2,1)$ and robot $\rob _1$ to $(r-\alpha+2,0)$. We split the next stage of movements into two cases depending on the parity of $r$. If $r$ is odd, perform the four moves $\rob _3: (2,1) \move (1,1)$, $\rob _1: (r-\alpha +2,0) \move (r-\alpha +1,0)$, $\rob _3: (1,1) \move (1,0)$, and $\rob _0: (r-\alpha ,1) \move (\alpha ,1)$ in that order. If $r$ is even perform the following sequence of five moves: $\rob _3: (2,1) \move (2,0)$, $\rob _1: (r-\alpha +2,0) \move (r-\alpha +1,0)$, $\rob _0: (r-\alpha ,1) \move (r-\alpha -1,1)$, $\rob _3: (2,0) \move (1,0)$ and $\rob _0: (r-\alpha -1,1) \move (\alpha ,1)$. In the final stage, perform the following three moves $\rob _1: (r-\alpha +1,0) \move (r-\alpha +1,1)$, $\rob _3: (1,0) \move (0,0)$ and $\rob _1: (r-\alpha +1,1) \move (r-\alpha ,1)$. It can be seen that all of these moves were legal. The effect of these three sequences of movements is to permute the positions of the robots $\rob _0,\rob _1, \rob _3$ in a 3-cycle and leave $\rob _2$ fixed. Considering the symmetry $(i,j) \rightarrow (-i,2-j)$ of the cylinder, we see that we can also permute the positions of the robots $\rob _0,\rob _2, \rob _3$ in a 3-cycle, leaving robot $\rob _1$ fixed. By combining these two 3-cycles we see that we can bring any of the four robots into any of the positions in the starting configuration, and so the set was a completely mobile general position set. The result for $s \geq 4$ now follows from Lemma~\ref{lem:shifting}, completing the proof.
\end{proof}

\section{Mobile general position for torus graphs}\label{sec:twocycle}

In this section, we shall completely classify the mobile general position numbers of torus graphs, settling an open question from \cite{klavzar-2025}. We begin with the following results which can be easily confirmed.  
\begin{claim} \label{c1}
Let $u,v,z \in V(C_s \Box C_t)$. Then $z \in I(u,v)$ if and only  if $z_x \in [u_x,v_x]_s$ and 
$z_y \in [u_y,v_y]_t$.
\end{claim}

\begin{claim} \label{c2}
Let $u,v,z \in [s]_0$ and $a$ be a non-negative integer. 
Then 
 $z \in [u,v]_s$ if and only if $z + a \in [u + a,v +a]_s$, 
 where additions are performed modulo $s$. 
\end{claim}

Claim~\ref{c2} as well as the fact that $C_s \Box C_t$ is 
a vertex-transitive graph is the basis for the following definitions.

Let $S$ and $S'$ be subsets of $V(C_s \Box C_t)$ and 
$a,b$ be non-negative integers. 
If $u\in S$, then $u + (a,b)$
(where additions in the first and second coordinates are 
modulo $s$ and $t$ respectively) is called the \textit{shift} of $u$ by $(a,b)$. 
If $|S|=|S'|$ and for every $u\in S$   
$u + (a,b) \in S'$, i.e.\ the shift of $u$ by $(a,b)$ belongs to $S'$,
then we say that $S$ and $S'$ are \textit{symmetric} with respect to $(a,b)$.
By Claim~\ref{c2} we obtain
\begin{claim} \label{c3}
Let $S$ and $S'$ be symmetric subsets of $V(C_s \Box C_t)$
with respect to $(a,b)$.  
If $u,v,z \in S$ and $u', v', z' \in S'$ are their 
corresponding shifts by $(a,b)$, then 
$z \in I(u,v)$ if and only if  $z' \in I(u',v')$.
\end{claim}

Next, we recall the general position number of a Cartesian product of two cycles.

\begin{theorem}[\cite{KorzeVesel}]\label{gp}
Let $s\ge t\ge 4$. Then
\[
\gp(C_s\cp C_t)=
\begin{cases}
5, & t=4,\\
6, & t\in\{5,6\}\text{ or }s,t\in\{8,10,12\},\\
7, & \text{otherwise}.
\end{cases}
\]
\end{theorem}

\subsection{The base configuration on $C_7\cp C_7$}

\begin{proposition}\label{prop:C7C7}
\[
\cmob(C_7\cp C_7)=\mob(C_7\cp C_7)=6.
\]
\end{proposition}

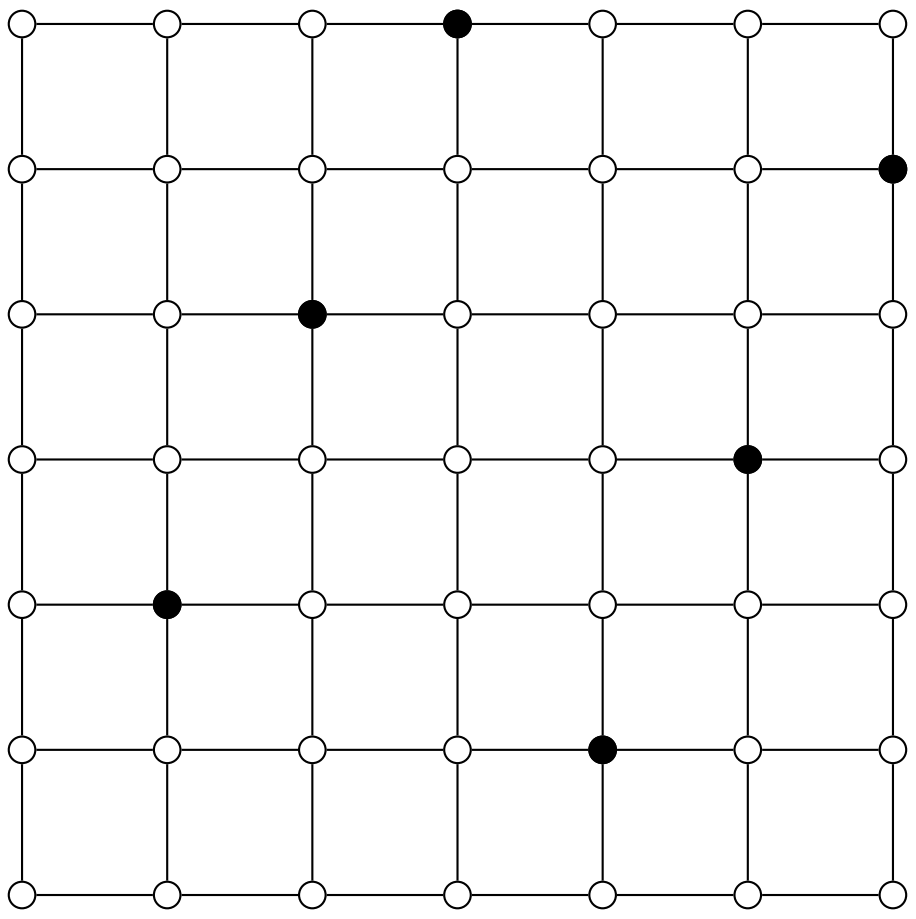
\begin{figure}[bt] 
\centering
\begin{tikzpicture}[scale=0.8]

 \tikzstyle{nod}= [circle, draw,inner sep=0pt, minimum size=5]
 
 \tikzstyle{nod2}= [circle, fill = black, draw,inner sep=0pt, minimum size=5]

  \foreach \x in {0,...,6}
    \foreach \y in {0,...,6} 
       {  
       \node [nod]  (\x\y) at (1.2*\x,1.2*\y) {};} 

  \foreach \x in {0,...,6}
    \foreach \y [count=\yi] in {0,...,5}  
      \draw (\x\y)--(\x\yi) (\y\x)--(\yi\x) ;

\foreach \pt in {12, 24, 36, 41, 53, 65} {
    \node[nod2] at (\pt) {};
  }
\end{tikzpicture}
\caption{A mobile general position set in $C_{7} \cp C_{7}$} \label{C7C7}
\end{figure}

\begin{proof}
We first prove the upper bound. Let $S$ be a general position set of cardinality $7$ in $C_7\cp C_7$. 
By Proposition \ref{prop:one-per-layer},  every layer contains exactly one vertex of $S$.
After any move of one robot, the resulting set would contain two robots in one layer and none in another; it therefore could not be a general position set of cardinality $7$. Consequently,
\[
\mob(C_7\cp C_7)\le6.
\]

For the reverse inequalities, consider
\[
S=\{(1,2),(2,4),(3,6),(4,1),(5,3),(6,5)\},
\]
shown in Fig.~\ref{C7C7}. A direct check using Claim~\ref{c1} shows that $S$ is in general position.

Move the robots one at a time by one step in the negative $x$-direction, in the order in which they are listed above. Every intermediate configuration is in general position, and after the six moves the entire configuration has been shifted by $(-1,0)$. Repeating this block of moves shifts the configuration through all seven positions in the $x$-direction. Similarly, we may move the robots in the same order in the negative $y$-direction. Every intermediate configuration is in general position and after six moves the configuration has been shifted by $(0,-1)$. Combining these two sequences of moves, we see that each robot may visit each vertex and $S$ is completely mobile. Therefore,
%
%
\[
\cmob(C_7\cp C_7)=\mob(C_7\cp C_7)=6.
\]
\end{proof}

\subsection{Stretched seven-vertex configurations}

For positive integers $a$ and $b$, let
\[
s^{a,b}=(s^{a,b}_0,\ldots,s^{a,b}_6)
:=\bigl((0,0),(a,2b),(2a,4b),(3a,6b),(4a,b),(5a,3b),(6a,5b)\bigr).
\]
For $k\in[7]_0$, write $s^{a,b}_k=(s^{a,b}_{k,x},s^{a,b}_{k,y})$. The entries of $s^{a,b}$ may be regarded as vertices of $C_{7a}\cp C_{7b}$. In particular,
\[
s^{1,1}=\bigl((0,0),(1,2),(2,4),(3,6),(4,1),(5,3),(6,5)\bigr)
\]
is a general position set in $C_7\cp C_7$.  Indeed, the distance between any two of its vertices belongs to $\{3,4,5\}$.

To extend this construction, define the following perturbation tuples in the $x$-direction:
\begin{align*}
{}^{1}\!\delta^{x}&=((0,0),(0,0),(0,0),(1,0),(1,0),(1,0),(1,0)),\\
{}^{2}\!\delta^{x}&=((0,0),(0,0),(0,0),(1,0),(2,0),(2,0),(2,0)),\\
{}^{3}\!\delta^{x}&=((0,0),(0,0),(1,0),(2,0),(3,0),(3,0),(3,0)),\\
{}^{4}\!\delta^{x}&=((0,0),(0,0),(1,0),(2,0),(3,0),(4,0),(4,0)),\\
{}^{5}\!\delta^{x}&=((0,0),(1,0),(2,0),(3,0),(4,0),(5,0),(5,0)),\\
{}^{6}\!\delta^{x}&=((0,0),(1,0),(2,0),(3,0),(4,0),(5,0),(6,0)),
\end{align*}
and the following perturbation tuples in the $y$-direction:
\begin{align*}
{}^{1}\!\delta^{y}&=((0,0),(0,0),(0,1),(0,1),(0,0),(0,1),(0,1)),\\
{}^{2}\!\delta^{y}&=((0,0),(0,0),(0,1),(0,2),(0,0),(0,1),(0,1)),\\
{}^{3}\!\delta^{y}&=((0,0),(0,1),(0,2),(0,3),(0,0),(0,2),(0,2)),\\
{}^{4}\!\delta^{y}&=((0,0),(0,1),(0,2),(0,4),(0,0),(0,2),(0,3)),\\
{}^{5}\!\delta^{y}&=((0,0),(0,2),(0,3),(0,5),(0,1),(0,3),(0,4)),\\
{}^{6}\!\delta^{y}&=((0,0),(0,2),(0,4),(0,6),(0,1),(0,3),(0,5)).
\end{align*}
For convenience, set ${}^{0}\!\delta^x={}^{0}\!\delta^y=((0,0),\ldots,(0,0))$. For $i,j\in[7]_0$, addition of the tuples is component-wise. We write
\[
S(a,b;i,j)
=\bigl\{s^{a,b}_k+{}^{i}\!\delta^x_k+{}^{j}\!\delta^y_k:k\in[7]_0\bigr\}.
\]
The first and second coordinates are interpreted modulo $7a+i$ and $7b+j$, respectively, so that
\[
S(a,b;i,j)\subseteq V(C_{7a+i}\cp C_{7b+j}).
\]
When $i=0$ or $j=0$, the corresponding perturbation is omitted from the notation.

It was proved in \cite{KorzeVesel} that $S(1,1;i,0)$ is a general position set in $C_{7+i}\cp C_7$ for every $i\in[7]$. We now show that these sets are mobile. For $i\in[7]$, let $m_x^{(i)}$ be the following permutation of $[7]_0$: 
\[ \begin{array}{c|c} i & m_x^{(i)}\\ \hline 
1&(2,1,0,6,5,4,3)\\ 
2&(2,1,0,6,5,4,3)\\ 
3&(1,0,6,5,2,3,4)\\ 
4&(1,0,6,5,4,3,2)\\ 
5&(4,3,2,1,0,6,5)\\ 
6&(0,1,2,3,4,5,6)\\ 
7&(0,1,2,3,4,5,6) \end{array} 
\] 
The entries specify the order in which the seven robots are moved one step in the positive $x$-direction. 
The index $i=7$ is used only for the case $C_{14}\cp C_7$. In this case the corresponding configuration is \[ S(2,1;0,0)=s^{2,1}. \]

\begin{figure}[bt] 
\centering
\begin{tikzpicture}[scale=0.7]
 \footnotesize 

 \tikzstyle{nod}= [circle, draw,inner sep=0pt, minimum size=5]
 
 \tikzstyle{nod2}= [circle, fill = black, draw,inner sep=0pt, minimum size=5]

  \foreach \x in {0,...,7}
    \foreach \y in {0,...,6} 
       {  
       \node [nod]  (\x\y) at (1.2*\x,1.2*\y) {};} 
\node[right=0.6 of 75] (85) {};

  \foreach \x in {0,...,7}
    \foreach \y [count=\yi] in {0,...,5}  
      \draw (\x\y)--(\x\yi);

  \foreach \y in {0,...,6}  
   \foreach \x [count=\xi] in {0,...,6}
      \draw (\x\y)--(\xi\y);

\foreach \pt in {00, 12, 24, 46, 51, 63, 75} {
    \node[nod2] at (\pt) {};
  }

\foreach \from/\to/\lbl in {
  46/56/7,
  51/61/6,
  63/73/5,
  75/85/4,
  00/10/3,
  12/22/2,
  24/34/1
} {
  \draw[->, >=latex, line width=1pt] (\from) -- (\to) node[midway, above]{\lbl};
}
  
\end{tikzpicture}
\caption{The first block of seven moves for a mobile general position set in $C_{8} \cp C_{7}$} \label{C8C7}
\end{figure}
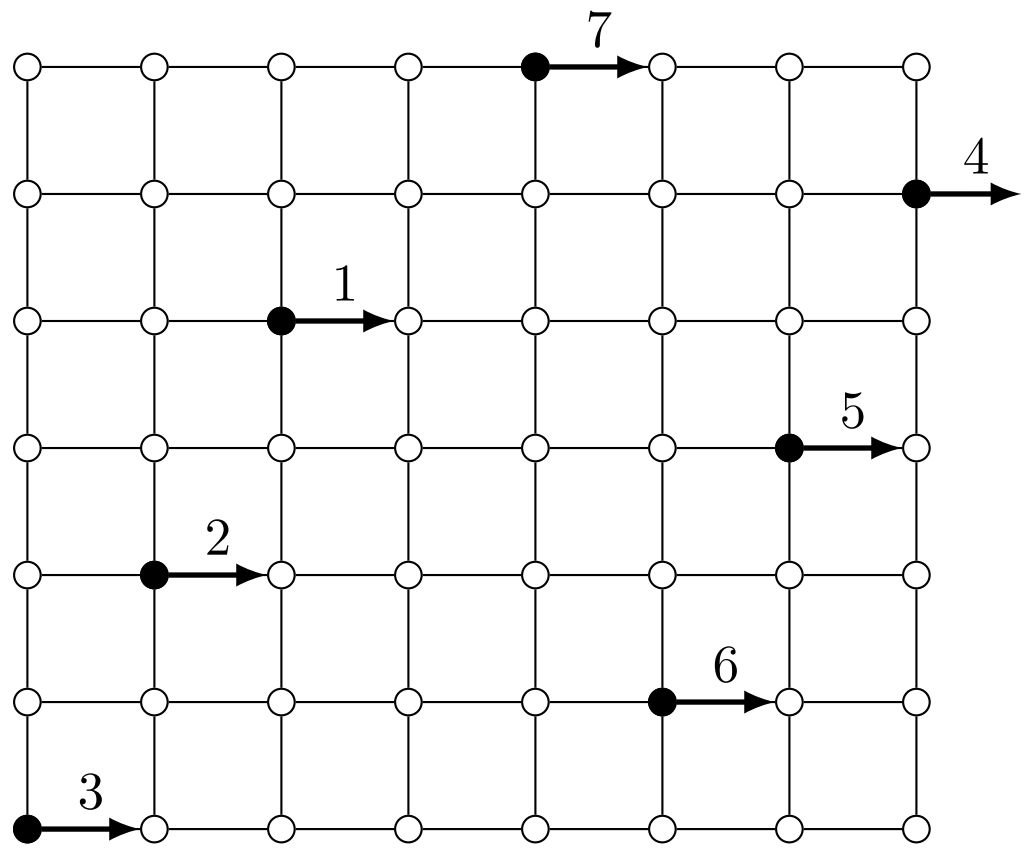

\begin{proposition}\label{prop:11} 
For every $s\in\{8,\ldots,14\}$ there exists a mobile general position set of cardinality $7$ in $C_s\cp C_7$. Consequently, 
\[ \mob(C_s\cp C_7)=7,\qquad 8\le s\le14. \] 
\end{proposition}
 \begin{proof} For $i\in[7]$, define 
\[ (a_i,r_i)= \begin{cases} (1,i),&1\le i\le6,\\ (2,0),&i=7, \end{cases} \] and put \[ S_i=S(a_i,1;r_i,0). \] 
Thus 
\[ S_i=S(1,1;i,0)\subseteq C_{7+i}\cp C_7 \qquad (1\le i\le6), \] 
while \[ S_7=S(2,1;0,0)=s^{2,1}\subseteq C_{14}\cp C_7. \] 
Starting from $S_i$, move the seven robots one at a time, each by one step in the positive $x$-direction, in the order $m_x^{(i)}$. 
A direct verification using Claim~\ref{c1} shows that every intermediate configuration is in general position. Figure~\ref{C8C7} illustrates the case $i=1$. After all seven moves, each robot has been shifted by $(1,0)$. Hence the resulting configuration is \[ S_i+(1,0), \] and the same block of seven moves can be repeated. Since the seven robots have pairwise distinct second coordinates, repeated horizontal sweeps visit every vertex of the corresponding graph $C_{7a_i+r_i}\cp C_7$. Therefore $S_i$ is a mobile general position set. It follows that 
\[ \mob(C_s\cp C_7)\ge7 \qquad\text{for }8\le s\le14. \] 
The reverse inequality follows from Theorem~\ref{gp}, and hence \[ \mob(C_s\cp C_7)=7. \] 
\end{proof}

The analogous construction in the $y$-direction requires the vertically stretched base tuple $s^{1,2}$. For $j\in[6]$, let
\[
\begin{array}{c|c}
 j & m_y^{(j)}\\ \hline
1&(1,4,0,3,6,2,5)\\
2&(1,6,4,2,0,5,3)\\
3&(4,6,2,0,5,1,3)\\
4&(1,6,4,2,0,5,3)\\
5&(4,6,2,0,5,1,3)\\
6&(1,6,4,2,0,5,3)
\end{array}
\]
be the order in which the robots are moved one step in the positive $y$-direction.

\begin{proposition}\label{prop:11y}
For every $j\in[6]$, the set $S(1,2;0,j)$ is a mobile general position set in $C_7\cp C_{14+j}$. Consequently,
\[
\mob(C_7\cp C_{14+j})=7.
\]
\end{proposition}

\begin{proof}
For each $j\in[6]$, a direct verification shows that the seven configurations obtained by moving the robots in the order $m_y^{(j)}$ are in general position. After the seventh move, the configuration has been shifted by $(0,1)$. Repeating the block therefore moves the configuration through all positions in the $y$-direction. Since its seven robots have pairwise distinct $x$-coordinates, every vertex of $C_7\cp C_{14+j}$ is visited. The upper bound again follows from Theorem~\ref{gp}.
\end{proof}

\subsection{Periodic move sequences}

We next formalise the preceding move sequences. For $i,j\in[6]$ and $p,q\in\{0,1,\ldots,7\}$, define
\[
\varepsilon^{i,p}_{x,k}=
\begin{cases}
1,&\text{if robot $k$ occurs among the first $p$ entries of $m_x^{(i)}$},\\
0,&\text{otherwise},
\end{cases}
\]
and define $\varepsilon^{j,q}_{y,k}$ analogously using $m_y^{(j)}$. The configuration after $p$ moves in the $x$-direction and $q$ moves in the $y$-direction is
\[
S^{p,q}(a,b;i,j)
=
\left\{
\left(
 a s^{1,1}_{k,x}+{}^{i}\!\delta^x_{k,x}+\varepsilon^{i,p}_{x,k},
 b s^{1,1}_{k,y}+{}^{j}\!\delta^y_{k,y}+\varepsilon^{j,q}_{y,k}
\right):k\in[7]_0
\right\},
\]
where the two coordinates are reduced modulo $7a+i$ and $7b+j$, respectively.

After a complete block of seven moves, the configuration is translated by one unit. Thus, if $p'=7u+p$ and $q'=7v+q$, where $p,q\in[7]_0$, then the configuration after $p'$ and $q'$ moves is a translate of $S^{p,q}(a,b;i,j)$ by $(u,v)$. By Claim~\ref{c3}, it is therefore enough to verify the $7\times7$ intermediate configurations with $p,q\in[7]_0$.

\begin{proposition}\label{prop:12}
Let $i\in[7]_0$ and $j\in[6]$. Then $S(1,2;i,j)$ is a mobile general position set in $C_{7+i}\cp C_{14+j}$.
\end{proposition}

\begin{proof} 
Fix $i\in[7]_0$ and $j\in[6]$. 

If $i\in[6]$, we use $m_x^{(i)}$ for the horizontal moves and $m_y^{(j)}$ for the vertical moves. Thus, starting from $S(1,2;i,j)$, a horizontal block consists of moving the seven robots, one at a time, one step in the positive $x$-direction in the order $m_x^{(i)}$, while a vertical block consists of moving them one at a time, one step in the positive $y$-direction in the order $m_y^{(j)}$. 

A direct computer verification, using Claim~\ref{c1}, shows that every intermediate configuration arising during these horizontal and vertical blocks is in general position. Equivalently, every configuration $S^{p,q}(1,2;i,j)$, with $p,q\in\{0,\ldots,7\}$, is in general position. 

After a complete horizontal block, all seven robots have moved once in the positive $x$-direction, so the configuration is translated by $(1,0)$. 
Similarly, after a complete vertical block it is translated by $(0,1)$. Hence the same two blocks can be repeated. By Claim~\ref{c3}, all translated configurations remain in general position. 
Since the translations $(1,0)$ and $(0,1)$ generate $C_{7+i}\cp C_{14+j}$, every vertex of the torus is eventually visited by a robot. 

It remains to consider $i=0$. In this case the initial configuration $S(1,2;0,j)$ already contains one robot in each of the seven $x$-coordinates, so no horizontal block is needed. We only use the vertical move order $m_y^{(j)}$. The computer verification shows that all intermediate configurations in a vertical block are in general position. 
After each complete block the configuration is translated by $(0,1)$. Repeating this block makes each of the seven robots visit all vertices in its fixed $x$-column. Since the seven robots occupy all seven $x$-coordinates, every vertex of $C_7\cp C_{14+j}$ is visited. 

Thus $S(1,2;i,j)$ is a mobile general position set in all cases. 
\end{proof}

For the scaling argument below, we use the following monotonicity lemma, adapted from \cite[Proposition~3.3]{KorzeVesel}.

\begin{definition}
Let $s$ be a positive integer. Two $s$-tuples
\[
(x_0,\ldots,x_{s-1})\quad\text{and}\quad(d_0,\ldots,d_{s-1})
\]
are \emph{stable} if the $x_i$ are pairwise distinct and
\[
x_i<x_j\quad\text{implies}\quad d_i\le d_j
\]
for all $i,j\in[s]_0$.
\end{definition}

\begin{proposition}\label{stable}
Let $s$ be odd, let $(x_0,\ldots,x_{s-1})$ and $(d_0,\ldots,d_{s-1})$ be stable, and suppose that $x_i,d_i,h\in[s]_0$ for every $i$. Let $f_A\ge f_B\ge s$. If
\[
f_Bx_i+d_i\notin[f_Bx_j+d_j,f_Bx_k+d_k]_{f_Bs+h},
\]
then
\[
f_Ax_i+d_i\notin[f_Ax_j+d_j,f_Ax_k+d_k]_{f_As+h}.
\]
\end{proposition}

For use in the general construction, set
\[
\Delta^{i,p}_{x,k}={}^{i}\!\delta^x_{k,x}
+\varepsilon^{i,p}_{x,k}-\varepsilon^{i,p}_{x,0},
\qquad
\Delta^{j,q}_{y,k}={}^{j}\!\delta^y_{k,y}
+\varepsilon^{j,q}_{y,k}-\varepsilon^{j,q}_{y,0}.
\]
Subtracting the same quantity from every coordinate only translates the configuration. Direct inspection shows that, for every $p, q$ , the tuples
\[
(s^{1,1}_{0,x},\ldots,s^{1,1}_{6,x})
\quad\text{and}\quad
(\Delta^{i,p}_{x,0},\ldots,\Delta^{i,p}_{x,6})
\]
are stable; the analogous statement holds for the $y$-coordinates. Moreover, all correction terms lie in $[7]_0$.

To include zero remainders in the scaling construction, use
\[
m_x^{(0)}=(6,5,4,3,2,1,0),
\qquad
m_y^{(0)}=(3,6,2,5,1,4,0).
\]
The first order is needed only when $a\ge2$, and the second only when $b\ge3$. 

For the graphs $C_{14}\cp C_{9}$, $C_{14}\cp C_{11}$ and $C_{14}\cp C_{13}$ we use 
the additional initial configurations and corresponding move orders listed in Table \ref{tab:14}.

\begin{table}[htbp]
  \centering
  \begin{tabular}{l|l|c|c}
     $t$ & Initial 7-tuple in $C_{14}\cp C_t$ & $m_x$ & $m_y$\\
    \hline
$9$ & $(0,0),(4,2),(6,5),(7,8),(10,1),(11,4),(12,6)$ &$\!(0,3,2,6,1,5,4)$ & $\!(1,6,4,2,0,5,3)$ \\
$11$ & $(0,0),(1,3),(3,7),(6,9),(7,2),(8,4),(9,8)$ & $\!(2,6,1,5,0,4,3)$ & $\!(0,5,3,1,6,4,2)$ \\
$13$ & $(0,0),(1,9),(2,6),(3,2),(5,11),(8,7),(10,3)$ & $\!(3,2,1,4,5,0,6)$ & $\!(4,5,0,2,6,1,3)$ \\
    \hline
  \end{tabular}
  \caption{Seven-robot configurations and corresponding move orders in some products with $C_{14}$}
\label{tab:14}
\end{table}

\begin{proposition}\label{general}
Let $a\ge1$, $b\ge2$, and $i,j\in[7]_0$, such that \((7a+i,7b+j)\not \in \{(8,14),(10,14), (12,14),(14,14)\}\). Then
\[
\mob(C_{7a+i}\cp C_{7b+j})=7.
\]
\end{proposition}

\begin{proof} 
The upper bound \[ \mob(C_{7a+i}\cp C_{7b+j})\le 7 \] follows from Theorem~\ref{gp}. We prove the reverse inequality by showing that the move sequences defined above can be performed while preserving general position. For $p,q\in[7]_0$, consider the configuration obtained after $p$ moves of the current horizontal block and $q$ moves of the current vertical block. 
Up to a translation, the coordinates of robot $k$ can be written as 
\[ X_k(a)=a s^{1,1}_{k,x}+\Delta^{i,p}_{x,k} \pmod{7a+i} \] and \[ Y_k(b)=b s^{1,1}_{k,y}+\Delta^{j,q}_{y,k} \pmod{7b+j}, \] where the correction terms $\Delta^{i,p}_{x,k}$ and $\Delta^{j,q}_{y,k}$ were defined above. Recall that the pairs of $7$-tuples 
\[ \bigl(s^{1,1}_{0,x},\ldots,s^{1,1}_{6,x}\bigr), \qquad \bigl(\Delta^{i,p}_{x,0},\ldots,\Delta^{i,p}_{x,6}\bigr) \] and 
\[ \bigl(s^{1,1}_{0,y},\ldots,s^{1,1}_{6,y}\bigr), \qquad \bigl(\Delta^{j,q}_{y,0},\ldots,\Delta^{j,q}_{y,6}\bigr) \] 
are stable. 
Hence Proposition~\ref{stable} can be applied independently to the two coordinates. The required general position property has been verified directly for the  range 
\[ 1\le a\le 7,\; 2\le b\le 7, \; i,j\in[7]_0, \; (7a+i,7b+j)\not \in \{(8,14),(10,14), (12,14),(14,14)\}, \] 
and for every intermediate configuration arising from the prescribed move orders. 
We now extend this finite verification to arbitrary $a\ge1$ and $b\ge2$. 
Put \[ \bar a=\min\{a,7\}, \qquad \bar b=\min\{b,7\}. \] 
Thus $(\bar a,\bar b)$ belongs to the finite range verified above. Fix $p,q\in[7]_0$ and three distinct robots $k_1,k_2,k_3$. Since the corresponding configuration for $(\bar a,\bar b)$ is in general position, Claim~\ref{c1} implies that at least one of the following two relations holds: 
\begin{equation}\label{eq:x-witness} X_{k_1}(\bar a)\notin [X_{k_2}(\bar a),X_{k_3}(\bar a)]_{7\bar a+i}, \end{equation} or 
\begin{equation}\label{eq:y-witness} Y_{k_1}(\bar b)\notin [Y_{k_2}(\bar b),Y_{k_3}(\bar b)]_{7\bar b+j}. \end{equation} 
Suppose first that \eqref{eq:x-witness} holds. If $a\le7$, then $\bar a=a$, so the same exclusion already holds in the desired cycle $C_{7a+i}$. If $a>7$, then $\bar a=7$, and Proposition~\ref{stable}, applied with \[ s=7,\qquad f_B=7,\qquad f_A=a,\qquad h=i, \] gives \[ X_{k_1}(a)\notin [X_{k_2}(a),X_{k_3}(a)]_{7a+i}. \] 
Hence the three robots are in general position in the target graph. The argument is identical if \eqref{eq:y-witness} holds. If $b\le7$, there is nothing to prove, since $\bar b=b$. If $b>7$, then $\bar b=7$, and Proposition~\ref{stable}, now with 
\[ s=7,\qquad f_B=7,\qquad f_A=b,\qquad h=j, \] 
yields 
\[ Y_{k_1}(b)\notin [Y_{k_2}(b),Y_{k_3}(b)]_{7b+j}. \] 
Again Claim~\ref{c1} shows that the three robots are in general position. 
Thus, for every ordered triple of distinct robots and for every intermediate configuration, at least one of the two coordinate exclusions required by Claim~\ref{c1} is preserved when $(\bar a,\bar b)$ is replaced by $(a,b)$. Consequently, every configuration occurring during the prescribed sequence of moves is a general position set. 
After one complete horizontal block, all seven robots have been translated by $(1,0)$, and after one complete vertical block they have been translated by $(0,1)$. 
Repeating these blocks therefore produces all translations of the configuration in \[ C_{7a+i}\cp C_{7b+j}. \] 
In particular, every vertex of the torus is eventually visited by a robot. Hence the initial set is a mobile general position set of cardinality $7$, and therefore \[ \mob(C_{7a+i}\cp C_{7b+j})\ge7. \] 
Together with the upper bound, this proves \[ \mob(C_{7a+i}\cp C_{7b+j})=7. \] \end{proof}

\subsection{Small torus graphs}
We now give the mobile general position numbers of torus graphs involving smaller cycles, which were not considered in the previous subsection.

\begin{table}[h]	
	\begin{tabular}{c|c|c|c}
		Torus & Initial Configuration & Move & Order \\ \hline  
		$C_3 \cp C_3$ & $(0,0), (0,1), (0,2)$ & $(i, j) \move (i+1, j)$ & $0,1,2$ \\
		$C_3 \cp C_4$ & $(0,0), (1,1), (0,2), (1,3)$ & $(i, j) \move (i-1, j)$ & $0,2,1,3$ \\
		$C_3 \cp C_5$ & $(0,0), (2,1), (0,2), (1,3), (2,4)$ & $(i, j) \move (i+1, j)$ & $0,4,3,2,1$ \\
		$C_3 \cp C_6$ & $(0,0), (1,1), (2,2), (1,4)$ & $(i, j) \move (i, j+1)$ & $2,3,1,0$ \\
		$C_3 \cp C_7$ & $(0,0), (1,1), (2,2), (0,3), (1,4), (2,5)$ & $(i, j) \move (i, j+1)$ & $5,4,3,2,1,0$ \\
		$C_4 \cp C_4$ & $(0,0), (2,1), (0,2), (2,3)$ & $(i, j) \move (i+1, j)$ & $0,1,2,3$ \\
		$C_4 \cp C_5$ & $(0,0), (3,1), (1,2), (0,3), (2,4)$ & $(i, j) \move (i+1, j)$ & $0,2,4,1,3$ \\
        $C_4 \cp C_6$ & $(0,0), (2,2), (0,3), (1,4), (3,4)$ & $(i, j) \move (i, j+1)$ & $4,3,0,2,1$ \\
		$C_5 \cp C_5$ & $(0,0), (1,1), (3,2), (4,3), (1,4)$ & $(i, j) \move (i+1, j)$ & $3,4,2,1,0$ \\
		$C_5 \cp C_6$ & $(0,0), (3,1), (4,2), (1,3), (2,4), (4,5)$ & $(i, j) \move (i+1, j)$ & $0,4,5,2,3,1$ \\
        $C_5 \cp C_7$ & $(0,0), (3,1), (4,2), (1,3), (2,4), (4,5)$ & $(i, j) \move (i, j+1)$ & $5,4,3,2,1,0$ \\
		$C_6 \cp C_6$ & $(0,0), (2,1), (1,2), (4,3), (3,4)$ & $(i, j) \move (i, j+1)$ & $2,1,4,3,0$ \\
        $C_6 \cp C_6$ & $(1,2), (2,1), (3,4), (4,3), (5,0)$ & $(i, j) \move (i, j+1)$ & $0,1,2,3,4$ \\
		$C_6 \cp C_7$ & $(0,0), (1,1), (4,2), (5,3), (2,4), (3,6)$ & $(i, j) \move (i, j+1)$ & $4,3,2,1,0,5$
	\end{tabular}	
	\caption{Initial configurations and moves for mobile general position sets of some small torus graphs. Guidance on using the table is given in the proof of Proposition~\ref{prop:c3tori}.}	
	\label{tab:smalltorusmoves}
\end{table}

\begin{proposition}
	For all integers $r \geq 3$, \[ \mob(C_3 \cp C_r) =  \begin{cases}
		3, & \text{if}\ r=3 \\
		4, & \text{if}\ r \in \{ 4, 6 \} \\
		5 & \text{if}\ r=5 \\
		6 & \text{if}\ r \geq 7.
	\end{cases} \]	
	\label{prop:c3tori}
\end{proposition}

\begin{proof}
	First, to see that $\mob(C_3 \cp C_3) \leq 3$, it suffices to note that if robots are stationed at $(0,0), (1,0)$, the only possible positions for other robots are $(2,1)$ and $(2,2)$, and no move may be performed from this configuration.
	
	To see that $\mob(C_3 \cp C_3) \geq 3$, station robots at the vertices and in the order indicated by Table~\ref{tab:smalltorusmoves}. Then repeatedly apply the move from the table to the robots in the order given. That is, move $(0,0) \move (1,0)$, $(0,1) \move (1,1)$, $(0,2) \move (1,2)$, $(1,0) \move (2,0)$, and so on. Eventually, all vertices will be traversed. It is easy to verify that general position is always maintained. 
	
	Similarly, applying Table~\ref{tab:smalltorusmoves} to $C_3 \cp C_4$ shows that $\mob(C_3 \cp C_4) \geq 4$. From~\cite[Theorem 3.4]{tian-2021a}, we know that $\gp(C_3 \cp C_4) = 4$, giving the upper bound. Table~\ref{tab:smalltorusmoves} also shows that $\mob(C_3 \cp C_5) \geq 5$. Since $\gp(C_3 \cp C_5) = 5$, we have equality.
	
	We now show that $\mob(C_3 \cp C_6) \leq 4$. Let $S$ be a general position set of $C_3 \cp C_6$ with $|S| = 5$. It is clear that $|S \cap C_3^i| \leq 1$ for $i \in [6]_0$. Also, if $(i,j), (i,k) \in S$ and $j \neq k$, $|j-k|_6 = 3$. Therefore, up to isomorphism, there is a single possibility for $S$. There are no legal moves, so we deduce the required result.
	


    For $r \geq 8$, let $S = \{ (0,0), (1,1), (2,2), (0, 3+ \lfloor\frac{r-6}{2}\rfloor), (1, 4 + \lfloor\frac{r-6}{2}\rfloor), (2, 5 + \lfloor\frac{r-6}{2}\rfloor) \}$ and move $(i, j) \move (i, j+1)$ in the reverse of the lexicographic order. It can be seen that at each stage the set is in general position, and it follows that $S$ is a mobile general position set. Again by~\cite[Theorem 3.4]{tian-2021a}, $\gp(C_3 \cp C_r) = 6$ and so we have $\mob(C_3 \cp C_r) = 6$.
\end{proof}

\begin{proposition}
	For all integers $r \geq 4$, \[ \mob(C_4 \cp C_r) =  \begin{cases}
		4, & \text{if}\ r=4, \\
		5, & \text{if}\ r > 4.
	\end{cases} \]	
    \label{prop:c4tori}
\end{proposition}
\begin{proof}
    First, we show that $\mob(C_4 \cp C_4) = 4$. The lower bound follows from Table~\ref{tab:smalltorusmoves}. To see that $\mob(C_4 \cp C_4) \leq 4$, let $S$ be a general position set of $C_4 \cp C_4$ with $|S| = 5$. Remembering that $\gp(C_4) = 2$, there must be a $C_4$-layer containing exactly $2$ robots and at most one $C_4$-layer that has no vertices in $S$. Suppose without loss of generality that $|C_4^0 \cap S| = 2$, $|C_4^1 \cap S| > 0$. It is easy to verify that we may assume $C_4^0 \cap S = \{ (0,0), (2,0) \}$. If $(i,j), (k,j) \in S$ with $i \neq k$, $S \cap ({}^iC_5 \cup {}^kC_5) = \{ (i,j), (k,j) \}$. Therefore, $S - \{(0,0), (2,0)\} \subseteq {}^1C_4 \cup {}^3C_4$ and for $k \in [3]$, $|C_4^k \cap S| = 1$. We can see that there is no possible legal move, giving the required result.

    For $4 < r \leq 6$, $\mob(C_4 \cp C_r) \leq \gp(C_4 \cp C_r) = 5$ \cite[Theorem 2.8]{KorzeVesel}. Table~\ref{tab:smalltorusmoves} gives the required lower bounds. 
    For $r \geq 7$, taking $S = \{ (0,0), (2,1), (0, 2 + \lceil\frac{r-5}{2}\rceil), (1, 3 + \lceil\frac{r-5}{2}\rceil), (3, 4 + \lceil\frac{r-5}{2}\rceil) \}$ and moving $(i, j) \move (i, j+1)$ in the reverse of the lexicographic order, we see that $S$ is a mobile general position set. Hence, $\mob(C_4 \cp C_r) = 5$.
\end{proof}

\begin{proposition}
	For all integers $r \geq 5$, \[ \mob(C_5 \cp C_r) =  \begin{cases}
		5, & \text{if}\ r=5, \\
		6, & \text{if}\ r \geq 6.
	\end{cases} \]	
\end{proposition}
\begin{proof}
    By Table~\ref{tab:smalltorusmoves}, $\mob(C_5 \cp C_5) \geq 5$, $\mob(C_5 \cp C_6)$ and $\mob(C_5 \cp C_7) \geq 6$.
    For $r \geq 8$ taking $S = \{ (0,0), (3,1), (4,2), (1, 3 + \lfloor\frac{r-6}{2}\rfloor), (2, 4 + \lfloor\frac{r-6}{2}\rfloor), (4, 5 + \lfloor\frac{r-6}{2}\rfloor) \}$ and moving $(i, j) \move (i, j+1)$ in the reverse of the lexicographic order, we see that $S$ is mobile general position set.
    Since $\gp(C_5 \cp C_r) = 6$ for $r \geq 6$~\cite[Theorem 2.9]{KorzeVesel}, it then follows immediately that $\mob(C_5 \cp C_r) = 6$.

    To see that $\mob(C_5 \cp C_5) \leq 5$, let $S$ be any general position set of $C_5 \cp C_5$ with $|S| = 6$. It is easy to see that each $|C_5^i \cap S| \leq 2$ for $i \in [5]_0$. If $(i,j), (k,j) \in S$ with $i \neq k$, $S \cap ({}^iC_5 \cup {}^kC_5) = \{ (i,j), (k,j) \}$. So, if there are $2$ pairs of vertices in $S$ with the same second coordinate, the remaining $2$ vertices must have the same first coordinate. However, one can verify that the vertices cannot be in general position. 
    
    Hence, we can assume without loss of generality that $|C_5^0 \cap S| = 2$ and $|C_5^i \cap S| = 1$ for $i>0$. Also, the vertices in $C_5^0$ cannot be adjacent, so we may assume that $(0,0), (0,2) \in S$. Label the vertex in $C_5^i \cap S$ as $(x_{i}, i)$ for $i \in [4]$. We must have either $x_1 = x_4 = 1$ or $x_1, x_4 \in \{ 3,4 \}$.

    If $x_1, x_4 \in \{ 3,4 \}$, we must have $x_2, x_3 = 1$. But $(x_1, 1), (1, 2), (1,3)$ is a shortest path, so we have a contradiction. Therefore, $x_1 = x_4 = 1$ and so $x_2, x_3 \in \{ 3,4 \}$. We cannot have $x_2 = x_3$, so this leaves two possibilities for $S$. However, there are no legal moves in either case, hence $\mob(C_5 \cp C_5) = 5$.
\end{proof}

\begin{proposition}
	For all integers $r \geq 6$, \[ \mob(C_6 \cp C_r) =  \begin{cases}
		5, & \text{if}\ r=6, \\
		6, & \text{if}\ r \geq 7.
	\end{cases} \]	
\end{proposition}
\begin{proof}
    Firstly, the lower bounds follow from the configurations in Table~\ref{tab:smalltorusmoves}. Notice that in the case of $C_6 \cp C_6$ there are two given configurations. We must perform the moves for the first configuration and then move $(0,0) \move (5,0)$ before carrying out the moves for the second in order to visit all vertices. For $r \geq 8$, we see that $S = \{ (3, 5 + \lceil \frac{r-6}{2} \rceil), (4, 2), (2, 4 +  \lceil \frac{r-6}{2} \rceil), (1, 1), (5, 3 + \lceil \frac{r-6}{2} \rceil), (0,0) \}$ is a general position set of $C_6 \cp C_r$. Also, moving $(i,j) \move (i,j+1)$ in the order written is a sequence of valid moves. Therefore, $S$ is a mobile general position set of $C_6 \cp C_r$ and $\mob(C_6 \cp C_r) \geq 6$.
    By~\cite[Theorem 2.11]{KorzeVesel} we have $\gp(C_6 \cp C_r) = 6$ for $r \geq 7$, so $\mob(C_6 \cp C_r) = 6$.

    Finally, we show that $\mob(C_6 \cp C_6) \leq 5$. We know from \cite[Theorem 2.11]{KorzeVesel} that $\gp(C_6 \cp C_6) = 6$. Let $S$ be a general position set of $C_6 \cp C_6$ with $|S| = 6$, $(0,0) \in S$. First we prove the following claim: if $(0,i) \in S$ for $i \in [5]$, then \[S = \{ (0,0), (k,1), (j,2), (0,3), (k,4), (j,5) \}\] where $\{ k,j \} = \{ 2,4 \}$.
    
    Suppose that $(0,i) \in S$ for some $i \in [5]$. By symmetry, we may assume that $i \in [3]$. Notice that $|C_6^0 \cap S| = |C_6^i \cap S| = 1$. 
    Since $(0,0), (0,1)$ is a maximal general position set of $C_6 \cp C_6$, $i \neq 1$.

    If $i=2$, $(C_6^3 \cup C_6^5) \cap S = \emptyset$. But then $C_6^1 \cap S = \{ (j_1,1), (k_1, 1) \}$ and $C_6^4 \cap S = \{ (j_2, 4), (k_2, 4) \} $ for some $j_1, j_2 \in \{ 1,2 \}$ and $k_1, k_2 \in \{ 4,5 \}$. One can verify that in any case, $S$ cannot be a general position set.

    Hence, $i=3$. Suppose $|C_6^1 \cap S| > 1$. Clearly, $C_6^1 \cap S = \{  (j,1), (k,1) \}$ where $j \in \{ 1,2 \}$, $k \in \{ 4,5 \}$ and $C_6^2 \cap S = \emptyset$. Therefore, $(C_6^4 \cup C_6^5) \cap S = \{ (x_1, y_1), (x_2, y_2) \}$. It is clear that $x_1, x_2 \neq 3$. Also, $x_1, x_2 \notin \{ 1,5 \}$. If $j=2$ and $k=4$, then $x_1, x_2 \notin \{ 2, 4 \}$. But if $j=1$, $x_1, x_2 \neq 4$ and if $k=5$, $x_1, x_2 \neq 2$. Therefore, we must have $x_1 = x_2$. But then $\{ (0,0), (x_1, y_1), (x_2, y_2) \}$ are not in general position --- a contradiction. Hence, we may write $S = \{ (x_i, i) : i\in[6] _0\}$, with $x_0 = x_3 = 0$. 
    
    It is clear that $x_i \neq 3$ for all $i \in [6]_0$. Suppose that $x_1 \in \{ 1, 2\}$. Then $x_2 \in \{ 4,5 \}$ and $x_4, x_5 \notin \{ 1,5 \}$. Since $x_4 \neq x_5$, $\{ x_4, x_5 \} = \{ 2,4 \}$. If $x_1=1$, $x_5 = 2$, giving $x_4 = 4$. But $(0,0)$ lies on a shortest $(4,4), (1,1)$-path, a contradiction. Therefore, $x_1 = 2$, giving $x_5 = 4$, $x_4 = 2$ and $x_2 = 4$. Our claim follows by symmetry.

    It follows that in any general position set $S$ of $C_6 \cp C_6$ with $|S|=6$, either $|^iC_6 \cap S| = |C_6^i \cap S| = 1$ for $i \in [6]_0$; or $|^iC_6 \cap S| = 1$ and $|C_6^i \cap S| \in \{ 0,2 \}$ for $i \in [6]_0$; or $|^iC_6 \cap S| \in \{ 0,2 \}$ and $|C_6^i \cap S| = 1$ for $i \in [6]_0$. Thus, there is no legal move from any such general position set and $\mob(C_6 \cp C_6) = 5$.
\end{proof}

\subsection{Medium-sized torus graphs}
We conclude this section by giving the mobile general position numbers of the remaining torus graphs, completing our classification for this family.

\begin{proposition}\label{prop:medium}
Let $14\ge s\ge t\ge8$. If $s,t\in\{8,10,12,14\}$, then
\[
\mob(C_s\cp C_t)=6;
\]
otherwise,
\[
\mob(C_s\cp C_t)=7.
\]
\end{proposition}

\begin{proof}
Suppose first that $s$ and $t$ do not both belong to $\{8,10,12,14\}$. For each pair $(s,t)$, Table~\ref{tab:medium} gives an initial configuration
\[
r=(r_0,\ldots,r_6)
\]
and two permutations $m_x,m_y$ of $[7]_0$.

For $p,q\in\{0,1,\ldots,7\}$, let $r^{p,q}$ be the configuration obtained by moving the first $p$ robots in the order $m_x$ one step in the positive $x$-direction and the first $q$ robots in the order $m_y$ one step in the positive $y$-direction. A computer verification confirms that every $r^{p,q}$ is in general position. After a full block in either direction, the configuration is translated by one unit.  Claim~\ref{c3} therefore shows that all subsequent configurations are also in general position. Repeating the blocks visits the entire torus and proves the lower bound $7$. Since the upper bounds follow from Theorem~\ref{gp}, the 
case for $s,t$ not both in $\{8,10,12,14\}$ is settled. 

\begin{table}[h]	
	\resizebox{\textwidth}{!}{%
\begin{tabular}{|c|c|c|c|c|}
	      \hline  
	  $s$ & $t$ & $r$ & $m_x$  & $m_y$ \\
      \hline 
        \hline  
        9 & 8 & $(0,0), (1,4),  (2,7), (4,2), (5,6), (6,1), (8,5)$ & 2, 5, 1, 4, 0, 3, 6  & 3, 5, 0, 2, 4, 6, 1 \\ 
      \hline  
    11 & 8 & $(0,0), (1,4), (3,1), (4,6), (6,2), (8,7),(9,3)$ & 3, 6, 2, 5, 1, 4, 0  & 1, 6, 4, 2, 0, 5, 3 \\ 
    \hline  
    13 & 8 & $(0,0), (2,5), (12,4), (6,6), (3,1), (9,7), (8,3)$ & 0, 3, 2, 4, 5, 1, 6 & 4, 0, 5, 3, 1, 2, 6 \\ 
    \hline  
    9 & 9 & $(0,0), (1,3), (3,7), (4,2), (5,6), (7,1), (8,4)$ & 1, 4, 0, 3, 6, 2, 5 & 2, 6, 4, 1, 3, 5, 0\\ 
    \hline  
    10 & 9 & $(0,0), (1,2), (3,5), (5,7), (6,1), (8,4), (9,6)$ & 4, 2, 1, 3, 0, 6, 5 & 3, 1, 6, 4, 2, 0, 5\\ 
    \hline  
    11 & 9 & $(0,0), (1,5), (2,8), (4,3), (6,7), (7,2), (8,6)$ & 6, 3, 2, 5, 1, 4, 0 & 0, 3, 2, 5, 4, 6, 1 \\ 
    \hline  
    12 & 9 & $(0,0), (4,2), (5,5), (7,8), (9,1), (10,4), (11,6)$ & 0, 3, 6, 2, 5, 1, 4 & 1, 6, 4, 2, 0, 5, 3 \\ 
    \hline  
    13 & 9 & $(0,0), (4,2), (6,5), (7,8), (10,1), (11,4), (12,6)$ & 0, 3, 6, 2, 5, 1, 4 & 1, 6, 4, 2, 0, 5, 3 \\ 
    \hline  
    11 & 10 &$(0,0), (1,6), (3,1), (5,8), (6,2), (8,9), (10,5)$ &1, 2, 4, 5,  0, 3, 6 & 4, 1, 2, 6, 0, 5, 3 \\ 
    \hline  
    13 & 10 &$(0,0), (2,5), (5,2), (6,7), (8,3), (9,8), (12,4)$ & 5, 1, 0, 3, 4, 6, 2 & 1, 0, 6, 5, 4, 2, 3 \\ 
    \hline  
    11 & 11 & $(0,0), (3,9), (5,5), (6,3), (8,10), (9,6), (10,4)$ & 0, 3, 6, 2, 5, 1, 4 & 5, 0, 2, 4, 6, 1, 3 \\ 
    \hline  
    12 & 11 & $(0,0),(1, 3), (2, 6), (3, 9), (7, 1), (8, 4), (11, 8)$ & 3, 5, 2, 4, 1, 0, 6 & 4, 2, 0, 5, 3, 1, 6 \\ 
    \hline  
    13 & 11 & $(0,0), (1,3), (2,7), (5,9), (6,2), (7,4), (8,8)$ & 6, 5, 4, 3, 2, 1, 0 & 5, 0, 3, 1, 6, 4, 2 \\ 
    \hline  
    13 & 12 &$(0,0), (1,6), (2,1), (4,8), (7,3), (8,9), (10,4)$ & 3, 2, 6, 5, 1, 4, 0 & 5, 3, 1, 6, 4, 2, 0 \\ 
    \hline  
    13 & 13 &$(0,0), (1,9), (2,6), (3,2), (5,11), (8,7), (10,3)$ & 3, 2, 6, 5, 1, 4, 0 & 6, 1, 3, 5, 0, 4, 2 \\ 
    \hline
    \end{tabular}%
}
	\caption{Initial configurations and move orders for the medium-sized torus graphs. The entries in the last two columns give the orders in which the robots are moved one step in the positive $x$- and $y$-directions, respectively.}	
	\label{tab:medium}
\end{table}

Now suppose that $s,t\in\{8,10,12,14\}$. Use the six-vertex configuration
\[
\begin{aligned}
r=\{&(0,0),(1,1),
(\lfloor s/3\rfloor,\lfloor2t/3\rfloor),
(\lfloor s/3\rfloor+1,\lfloor2t/3\rfloor+1),\\
&(\lfloor2s/3\rfloor,\lfloor t/3\rfloor),
(\lfloor2s/3\rfloor+1,\lfloor t/3\rfloor+1)\}.
\end{aligned}
\]
Move the robots in the orders
\[
m_x=(3,2,5,4,1,0),
\qquad
m_y=(5,4,3,2,1,0).
\]
A direct verification shows that all $6\times6$ intermediate configurations are in general position. Repeating the move blocks visits every vertex, and hence the set is mobile and 
the lower bound is 6. 

The upper bound for $$(s,t) \notin \{(14,8),(14,10),(14,12),(14,14)\}$$ follows from Theorem~\ref{gp}. For the remaining pairs, namely $$(s,t)\in \{(14,8),(14,10),(14,12),(14,14)\},$$ we performed an exhaustive computer search, as follows.

We first generated all seven-vertex general position sets containing $(0,0)$ in the graphs under consideration. We found that the graphs $C_{14}\cp C_8$, $C_{14}\cp C_{10}$, $C_{14}\cp C_{12}$, and $C_{14}\cp C_{14}$ admit exactly 14, 126, 602, and 4120 such sets, respectively.

We then verified that, in each of the graphs $C_{14}\cp C_{8}$, $C_{14}\cp C_{10}$, and $C_{14}\cp C_{12}$, all of these sets admit only moves that change the second coordinate. Moreover, in $C_{14}\cp C_{14}$, each such set admits moves either in the $x$-direction or in the $y$-direction, but never in both directions.

Consequently, none of the graphs $C_s \cp C_t$, for $(s,t)\in \{(14,8),(14,10),(14,12),(14,14)\}$, admits a mobile general position set of cardinality 7. Hence, the mobile general position number of each of these graphs is 6.
\end{proof}

The preceding results, together with the small cases established in the previous subsection, lead to the desired classification.

\begin{theorem}\label{thm:all-tori}
Let $s\ge t\ge3$. Then
\[
\mob(C_s\cp C_t)=
\begin{cases}
3, & s=t=3,\\
4, & t=3,\ s\in\{4,6\},\text{ or }s=t=4,\\
5, & t=4,\ s\ge5,\text{ or }t\in\{5,6\},\ s=t,\text{ or }t=3,s=5,\\
6, & t=3,\ s\ge7,\text{ or }t=5,\ s\ge6,\text{ or }s=t=7,\\
   & \hspace{1.7em}\text{or }s,t\in\{8,10,12,14\}, \text{ or } t =6,\ s\ge7\\
7, & \text{otherwise}.
\end{cases}
\]
\end{theorem}

\section*{Acknowledgments}
During the preparation of this manuscript, we learned that Wenlin Zhang and Haobo Ma of the Omega Institute have independently determined the mobile general position numbers of the cylinders $C_7 \cp P_s$ and $C_{10} \cp P_s$ using AI.

The third author was supported by the Slovenian Research Agency under the grants P1-0297,  J1-70016,  J1-70045 and N1-0431.

\end{document}